# A parameterization of the Fermat curves satisfying $x^{2N} + y^{2N} = 1$

Kerry M. Soileau
July 19, 2007

**Abstract**

Note that the family of closed curves $C_N = \{(x, y) \in \mathbb{R}^2; x^{2N} + y^{2N} = 1\}$ for $N = 1, 2, 3, \cdots$ approaches the boundary of $[-1,1]^2$ as $N \to \infty$. In this paper we exhibit a natural parameterization of these curves and generalize to a larger class of equations.

## Discussion

Note that the family of closed curves $C_N = \{(x, y) \in \mathbb{R}^2; x^{2N} + y^{2N} = 1\}$ for $N = 1, 2, 3, \cdots$ approaches the boundary of $[-1,1]^2$ as $N \to \infty$. The following figure shows $C_1$ through $C_5$; $C_1$ is the unit circle and $C_5$ is outermost.

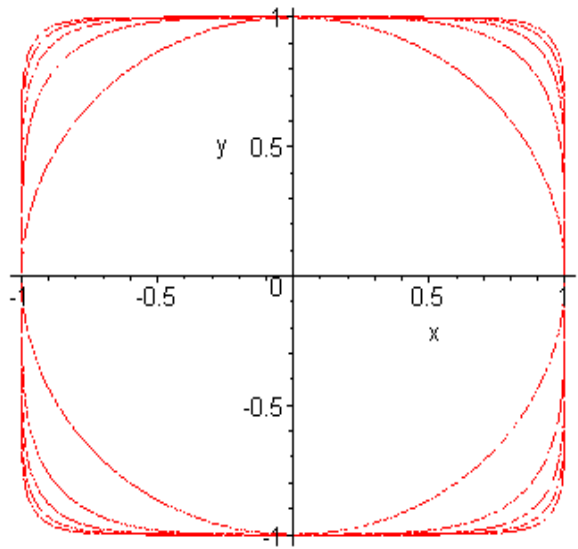

Let us denote the boundary of $[-1,1]^2$ by $C_\infty$. There is a natural parameterization of each curve $C_N$:

Kerry M. Soileau

$$x_N(\theta) = \left(\cos^{2N}\theta + \sin^{2N}\theta\right)^{-\frac{1}{2N}} \cos\theta$$

$$y_N(\theta) = \left(\cos^{2N}\theta + \sin^{2N}\theta\right)^{-\frac{1}{2N}} \sin\theta$$

for $0 \leq \theta < 2\pi$.

We may generalize this result to families of curves defined by the equation

$$(\alpha x + \beta y + \gamma)^{2N} + (\delta x + \varepsilon y + \zeta)^{2N} = 1$$

Let $L_N = \left\{(x,y) \in \mathbb{R}^2 ; (\alpha x + \beta y + \gamma)^{2N} + (\delta x + \varepsilon y + \zeta)^{2N} = 1\right\}$

If $\alpha\varepsilon \neq \beta\delta$, then

$$x_N(\theta) = \frac{1}{\alpha\varepsilon - \beta\delta}\left(\left(\cos^{2N}\theta + \sin^{2N}\theta\right)^{-\frac{1}{2N}}(\varepsilon\cos\theta - \beta\sin\theta) + \beta\zeta - \gamma\varepsilon\right)$$

$$y_N(\theta) = \frac{1}{\alpha\varepsilon - \beta\delta}\left(\left(\cos^{2N}\theta + \sin^{2N}\theta\right)^{-\frac{1}{2N}}(\alpha\sin\theta - \delta\cos\theta) + \delta\gamma - \alpha\zeta\right)$$

As $N \to \infty$ the $L_N$ approach a linear transformation of $C_\infty$; this mapping takes the point $(x,y)$ in $C_\infty$ to the point $\left(\dfrac{-\alpha\zeta + \gamma\delta - \delta x + \alpha y}{\alpha\varepsilon - \beta\delta}, \dfrac{\beta\zeta - \varepsilon\gamma + \varepsilon x - \beta y}{\alpha\varepsilon - \beta\delta}\right)$.

International Space Station Program Office, Avionics and Software Office, Mail Code OD, NASA Johnson Space Center, Houston, TX 77058
E-mail address: ksoileau@yahoo.com